\documentclass[journal]{IEEEtran}
\usepackage{graphics}

\font\MyScript=eusm10

\begin{document}
\title{An Algorithm for Optimal Partitioning of\\ 
Data on an Interval}

\author{Brad~Jackson,
Jeffrey~D.~Scargle,
David~Barnes, Sundararajan~Arabhi, Alina~Alt, 
Peter~Gioumousis, Elyus~Gwin, Paungkaew~Sangtrakulcharoen, 
Linda~Tan, and Tun~Tao~Tsai}

\maketitle

\thanks
{Manuscript received XXX, 2003.
This work was supported by the NASA Applied Information
Systems Research Program, the Woodward Fund of San Jose State University,
and the NASA Faculty Fellowship Program at Ames Research Center and Dryden Flight Research Center.}

\thanks
{B. Jackson is  with the Department~of~Mathematics, San~Jose~State~University,
J. Scargle is with the Space~Science~Division, NASA~Ames~Research~Center.  The other authors were participants
in the Center for Applied Mathematics and
Computer Science program at San Jose State University.}


\begin{abstract}
Many signal processing problems can be solved by 
maximizing the fitness of a segmented model over
all possible partitions of the data interval.
This letter describes a simple but powerful 
algorithm that searches
the exponentially large space 
of partitions of $N$ data points
in time $O(N^2)$.
The algorithm is guaranteed to 
find the exact global optimum, automatically
determines the model order (the number of segments),
has a convenient real-time mode, can be extended
to higher dimensional data spaces, and
solves a surprising variety of problems in signal
detection and characterization, density estimation,
cluster analysis and classification.
\end{abstract}


\noindent
EDICS: 1.STAT

\section{Introduction: The Problem}
A variety of signal processing and related problems 
can be viewed as the search 
for an optimal partition of data 
given on a time interval $I$.
For example, one may estimate a 
segmented model 
by maximizing some measure
of model fitness\footnote{This 
concept goes by many names, 
including \emph{fitness}, 
\emph{cost},
\emph{goodness of fit}, 
\emph{loss}, \emph{penalty}, 
\emph{objective function}, 
\emph{risk} {\it etc.}
Here we assume that 
the data analysis problem 
can be phrased in terms of 
some fitness that is to be maximized.}
defined on partitions of $I$.
Since the space of all partitions 
of a continuum is infinite,
it is advantageous to 
discretize the interval.
Often the data points themselves, say
\begin{equation}
X_{n}, n = 1, 2, \dots, N,
\end{equation}
\noindent
naturally subdivide $I$ into 
subintervals -- which we call {\it data cells}.
We avoid a precise definition 
because many types of data cells are possible.
Common examples are counts in bins,
measurements at a set of sample times
(evenly spaced or not), 
and event or point data.
The underlying idea is that 
restricting consideration to the finite space
of partitions whose elements are sets of
data cells will result in 
no significant loss 
of information or 
of resolution in the independent variable.

A {\it partition} ${\bf P}$ of an interval $I$ 
is a set of $M$ {\it blocks}, 
\begin{equation}
{\bf P}( I ) = \{ B_{m}, m \in \mbox{\MyScript{M}} \} \ ,
\mbox{\MyScript{M}} \equiv \{ 1, 2, \dots M \} \ , 
\label{blocks}
\end{equation}
\noindent
where the blocks are sets of data cells defined by
index sets $\mbox{\MyScript{N}}_{m}$:
\begin{equation}
B_{m} =  \{ X_{n}, n \in \mbox{\MyScript{N}}_{m} \}
\label{cells}
\end{equation}
\noindent
satisfying the usual conditions, $\bigcup\limits_{m} B_{m} = I$
and $B_{m} \bigcap B_{m'} = \emptyset$ if $m \ne m'$.
Here $\emptyset$ means a set of zero area,
since the 1D boundaries are not
relevant here.
We impose connectedness too---i.e.
that there be no gaps between the cells
comprising a given block.
$M$, the number of blocks,
must satisfy $0 \le M \le N$.
Partitions will be denoted in boldface, and 
refer to the interval $I$ unless otherwise stated.
Define
$P^{*}$ to be the (finite) set of all
partitions of $I$ into blocks.

Take as given an 
additive {\it fitness function}
that assigns
a value to any partition ${\bf P} \in P^{*}$ in the form
\begin{equation}
V({\bf P}) = \sum_{m=1}^{M} g(B_m) \ ,
\label{addprop}
\end{equation}
\noindent 
where $g(B_{m})$ is the fitness of block $B_{m}$.
Computationally, the data cells must
be represented by a data structure that contains 
{\it sufficient statistics} for the model -- {\it i.e.}
all information necessary to determine $g$
for any block.

We exhibit an efficient $O(N^2)$ dynamic programming algorithm 
that finds an {\it optimal partition}
${\bf P}^{\mbox{\small max}} \in P^{*}$:
$V({\bf P}^{\mbox{\small max}}) \ge V({\bf P})$ 
for all partitions ${\bf P} \in P^{*}$.  



Scargle \cite{scargle} proposed two greedy iterative 
algorithms for finding near-optimal partitions: 
one top-down (optimally divide $I$ into two parts; 
recursively do the same to each such part)
the other bottom-up (merge adjacent data cells).
In both cases Bayesian model comparison 
provides effective fitness functions and halting criteria,
implementing an $O(N^2)$ procedure for data spaces 
of 1, 2 and higher dimensions---hence the
term {\it Bayesian Blocks} \cite{scargle}.
But in practice these greedy algorithms 
often find significantly suboptimal partitions,
motivating the development reported here.

The dynamic programming idea 
for this kind of problem 
seems to have originated with
Richard Bellman \cite{bellman_2},
a paper which influenced Kay's
work \cite{kay_1,kay_2,kay_3}.
An extensive discussion of 
precisely the same problem addressed here,
but with a different approach to its solution,
is in \cite{barry_hartigan_1,barry_hartigan_2}.
Work by Hubert \cite{hubert_1,hubert_2},
with applications to meteorology,
influenced Kehagias and co-workers
\cite{fragkou,kehagias_1,kehagias_2,kehagias_3,kehagias_4,kehagias_5}, 
who developed a dynamic programming algorithm
much like ours, for applications such
as text segmentation (see also \cite{heinonen}), where
the raw data are provided in the form of
a similarity matrix.
\cite{Vidal} gives an $O(kN^2)$ dynamic 
programming algorithm for finding the 
optimal partition of an interval 
into $k$ blocks, for a given $k$.
See also \cite{igles,auger} 
for related work.
Thus while similar algorithms have been
previously developed,
ours finds the exact global optimum 
for any block-additive fitness function, 
automatically determines the number of segments, 
and can be used in either real-time 
or retrospective analysis.

\section{Dynamic Programming: Finding Optimal Partitions}

We describe an $O(N^2)$ 
algorithm that is guaranteed to 
solve the above problem by finding 
an exact global optimum,
for any fitness function $V$ that is additive 
in the sense of Eq. (\ref{addprop}).    
There is a large ($2^{N-1}$) but finite number 
of partitions in $P^{*}$.  
Dynamic programming \cite{bellman_1} 
is an intelligent method of searching 
this space of all possible solutions.  
Our algorithm can be applied whenever any subpartition 
of an optimal partition is optimal.
\newtheorem{theorem}{Theorem}
\begin{theorem}[Principle of Optimality] 
Let ${\bf P}^{\mbox{\small max}}$ be an 
optimal partition of $I$ and 
${\bf P}_{1} = \{ B_{m}, m \in a \}$ 
be any subset of the blocks of ${\bf P}^{\mbox{\small max}}$.  
Then ${\bf P}_1$ is an optimal partition of the
part of $I$ it covers, namely
$I_1 = \bigcup\limits_{m \in a} B_{m}$. 
\end{theorem} 

Intuitively this result  
follows from the contradiction that
a better subpartition of $I_{1}$
could be used to construct a partition of 
$I$ better than ${\bf P}^{\mbox{\small max}}$. 
The proof relies on the fact that the
block-additivity of the fitness function 
implies that it is also additive 
on subpartitions.  
To see this, divide partition ${\bf P}$ 
into any two disjoint parts,
${\bf P}_{1} = \{ B_{m} , m \in a\}$ 
and
${\bf P}_{2} = \{ B_{m} , m \in b\}$,
with
${\bf P}_{1} \bigcup {\bf P}_{2} = {\bf P}$ and 
$a \bigcup b = \mbox{\MyScript{M}}$.
Then the additivity of $V$ yields
\begin{eqnarray}
V({\bf P})  & = & \sum_{m=1}^{M} g(B_m) \nonumber \\
        & = & \sum_{m\in a} g(B_m) + \sum_{m\in b} g(B_m) \nonumber \\
        & = & V({\bf P}_1) + V({\bf P}_2). 
\end{eqnarray}

\noindent
{\bf Proof 1}: 
As above, denote by ${\bf P}_2$ the subpartition 
of ${\bf P}^{\mbox{\small max}}$, 
consisting of the blocks 
$\{B_m, m \in \mbox{\MyScript{M}} - a\}$
in ${\bf P}^{\mbox{\small max}}$ that are not in ${\bf P}_1$.
Let ${\bf P}_3$ be any other partition of $I_1$.  
Since ${\bf P}^{\mbox{\small max}}$ 
is an optimal partition of $I$ and ${\bf P}_3 \cup {\bf P}_2$ is also a partition 
of $I$ it follows that
$V({\bf P}^{\mbox{\small max}}) = V({\bf P}_1) + V({\bf P}_2) \ge V({\bf P}_3 \cup {\bf P}_2) = V({\bf P}_3) + V({\bf P}_2)$
so $V({\bf P}_1) \ge V({\bf P}_3).$  Thus ${\bf P}_1$ is an optimal partition of $I_1$.  

Dynamic programming is a recursive procedure that can be used to 
efficiently find the solution to many kinds of combinatorial 
optimization problems.  
Our algorithm derives the optimal partition of the first $n+1$ 
data points using previously obtained optimal partitions,
\emph{i.e.} those of the first 
$1, 2, \dots n$ data points.
At each iteration
we must consider all possible 
starting locations $j$, $1 \le j \le n$
of the last block  of the optimal partition.
For each putative $j$ the fitness function is
-- by the principle of optimality -- the
fitness of the optimal subpartition prior to $j$
plus the fitness of the last block itself.
The former was stored at previous iterations,
and the latter is a simple evaluation of $V$.
The desired new optimal partition corresponds
to the maximum over all $j$.

More precisely, 
define {\it opt}$(n)$  to be 
the value of the fitness function of the 
optimal partition ${\bf P}^{\mbox{\small max}}_{n}$ of 
the first $n$ cells of $I$, for $1 \le n \le N$.
The following dynamic programming algorithm 
finds the optimal partition ${\bf P}^{\mbox{\small max}}_{N}$: 
\begin{enumerate}
\item Define {\it opt}$(0) = 0$

\item Given that {\it opt}$(j)$ has been determined for $j = 0, 1, \dots , n$:
   \begin{itemize}
   \item Define {\it end}$(j,n+1) =  g(B_{j,n+1})$;
   $B_{j,n+1}$ is the union of cells $j, j+1, \dots , n+1$ 

   \item  Then compute
   \begin{equation}
   {\it opt}(n+1) = 
   \mathop{\mbox{Max}}_{j}\{opt(j-1) + \nonumber \\ {\it end}(j,n+1) \}, 
   \label{desired}
   \end{equation}
   \noindent
   for $j = 1, 2, \dots, n+1$.  

   \item
   The value of $j$ where this maximum occurs 
   is stored as {\it lastchange}$(n+1)$. 

   \end{itemize}

\item Repeat 2 until $n + 1 = N$, when
{\it opt}$(N)$, 
the optimal partition fitness for all $N$ cells, 
has been obtained.   

\item Backtrack using the {\it lastchange} vector 
to identify the start points of individual  
blocks of the optimal partition ${\bf P}^{\mbox{\small max}}$
in the following way.  
Let $n_1 = {\it lastchange}(N), n_2 = 
{\it lastchange}(n_1-1)$, {\it etc}. 
Then the last 
block in 
${\bf P}^{\mbox{\small max}}$ contains cells $n_1, n_{1}+1, \dots , N$, 
the next-to-last block 
in ${\bf P}^{\mbox{\small max}}$ contains 
cells $n_2, n_2 + 1, \dots , n_1 - 1,$ and so on.  
\end{enumerate}

\newtheorem{theorem_2}{Theorem}
\begin{theorem} This deterministic $O(N^2)$ dynamic programming 
algorithm finds the partition of $I$ that maximizes 
the (additive) fitness function.
\end{theorem}  

\noindent
{\bf Proof 2:} 
The proof is by mathematical induction.
Clearly 
${\it opt}(1) = 
{\mbox{Max} \{ {\it 0 + end}(1,1)} \} = g(B_{1,1})$
is the fitness of the only possible (and therefore optimal)
partition of the set comprising the first cell.  
At iteration $n+1$, 
assume not only that we have found the optimum
partition of ${\bf P}^{\mbox{\small max}}_{n}$, 
but also that for $i = 1, 2, \dots, n,$
we have stored the corresponding fitness
for this and all previous iterations in the 
array ${\it opt}(i)$,
and the index of the cell beginning this partition's
last block in array ${\it lastchange}(i)$.
Let $F(j) = opt(j-1) + {\it end}(j,n+1)$;
then the principle of optimality
shows that 
when $j$ indexes the first cell of the
last block of the desired partition 
${\bf P}^{\mbox{\small max}}_{n+1}, F(j)$ 
is the corresponding maximum fitness.
Further, for any $j, F(j)$
is the fitness of a legitimate partition of 
${\bf P}^{\mbox{\small max}}_{n+1}$,
namely that consisting of the optimal partition of the 
cells prior to $j$ followed by the single block 
$B_{j,n+1}$.
These two facts combine to prove that 
the maximum of $F(j)$ specified in Eq. (\ref{desired})
gives the desired optimum partition 
at iteration $n+1$.
Identification of the 
corresponding optimal blocks -- starting with the
last one and working backwards, as in 
part (4) of the algorithm -- can be 
validated with straightforward recursive 
application of the principle of 
optimality.  Finally, since 
{\it end}$(j,n+1) = g(B_{j,n+1})$ the algorithm 
requires $1 + 2+ \dots + N = O(N^2)$ evaluations 
of the function $g$.  
It also requires $O(N^2)$ additions and $O(N^2)$ 
comparisons in determining 
the maximums.

\section{Applications}

These results apply to any 
segmented modeling of 1-dimensional data
defined by a fitness function
that satisfies Eq. (\ref{addprop}).
Piecewise constant, or step functions 
form the most natural model class.
However the nature of the model
depends on the application, and many
other forms are possible, including
piecewise linear and piecewise 
exponential. 
The key is that the fitness function
must not depend on any model parameters
other than the changepoint locations.
We have found excellent results with
the posterior probability of the 
model for each segment,
given the data in that segment, 
marginalized over all parameters but
these locations \cite{scargle}.
Using logs turns the product posterior 
(resulting from statistical independence 
of the blocks)
into the required additive 
form Eq.(\ref{addprop}).

A comment about smoothness constraints
is in order. 
With some fitness functions 
the algorithm produces the 
degenerate solution assigning 
a segment to each data point.
In the Bayesian setting described above,
a natural way to address this problem
is to adopt a
prior distribution for the number of
segments, for example giving higher
weight to smaller numbers.
Indeed, the geometric prior \cite{coram}
corresponds to a constant term in the
fitness function for each block,
so there is no problem with maintaining
additivity.
This artifice controls model complexity,
but without imposing an explicit smoothness
condition with 
concomitant loss of time resolution.
Time series features on any time scale,
no matter how short, 
can be found if the data support them
in a statistically significant way.

Finally, we mention a few sample applications.
Implementing density estimation with 
piecewise constant Poisson models
yields histograms in which the bins
are not constrained to be equal.
The number of bins and their 
sizes and locations
are determined by the data.
The same model provides denoising and
structure estimation for time series
of events or counts of events in bins \cite{scargle}.

Further, almost all of the results described here
can be easily extended -- almost without change --
to data of higher dimensionality, as will
be described in future papers.
In this setting cluster analysis 
can be effected as a post-processing of
segmented models -- piecing the blocks
together into clusters -- and similarly with 
unsupervised classification and 
other data mining procedures.

\section{Conclusion}

As we have seen dynamic programming 
gives a good (polynomial) algorithm 
for finding an optimal partition of data on an 
interval for any fitness 
function $V$ satisfying the additive property [see Eq.(\ref{addprop})].
Ironically it has the same $O(N^2)$ complexity 
as the greedy algorithm.  
               
In comparing the use of our algorithm
to detect and characterize clusters 
(collections of blocks) with some of the standard 
clustering techniques \cite{Alpert}, 
we note that our method inherently compares 
partitions that have different numbers of blocks, 
so the number of blocks 
is automatically determined by the data. 
This is to be contrasted with most standard 
clustering techniques, in which $k$, 
the fixed number of clusters must be 
specified ahead of time.  
One often seeks to minimize the maximum 
diameter (defined as the maximum distance 
between any pair of points in 
the cluster) of the clusters, or to maximize 
the minimum separation between 
the clusters.  In dimension 1, 
there are well-known $O(kN^2)$ dynamic 
programming algorithms for finding the best 
partitions into k clusters.  
For dimension 2 and higher it is known 
that these standard problems 
are NP-complete.  We don't yet know if 
our problem is NP-complete in 
dimension 2 and higher.

In addition, considered as 
a density estimation or signal detection technique, 
our approach does not introduce any explicit smoothing
of the data.  Structure on any time scale, no
matter how short, will be detected if it is supported
by the data. 
While the parameter in the geometric prior
discussed above controls to some extent the
number of blocks -- and thus affects the roughness
of the optimized model -- it is not explicitly a
smoothing parameter.
Another feature is that the 
incremental way the algorithm operates on the
data makes a real-time mode trivial to implement.
This mode has found to be very useful in the
rapid detection of change points in a data stream.
And since {\it opt}$(i+1)$ is calculated from 
{\it opt}$(j), j = 1, 2, \dots , i,$ some 
of the necessary calculations can be performed 
as the data are still being collected.  
Also it is easy to modify the dynamic programming 
to yield the optimal partition with blocks of a 
minimum size 
(each block contains at least d data points, 
for a given positive integer d).  
These and other features are described in more detail at
an algorithm repository at:\\
\noindent
{\tiny
\verb+http://astrophysics.arc.nasa.gov/~pgazis/CodeArchiveServer.html+}

We are grateful to Steve Kay,
Thanasis Kehagias,
and the referees 
for helpful comments
and pointers to relevant 
earlier work.

\end{document}